\documentclass{amsart}

\newtheorem{theorem}{Theorem}[section]
\newtheorem{lemma}[theorem]{Lemma}

\newtheorem{bbcor}[theorem]{Corollary}
\newtheorem{proposition}[theorem]{Proposition}

\theoremstyle{definition}
\newtheorem{definition}[theorem]{Definition}

\theoremstyle{remark}
\newtheorem{remark}[theorem]{Remark}

\numberwithin{equation}{section}

\usepackage{graphicx}

\begin{document}

\title{Asymptotically maximal families of hypersurfaces in toric varieties}

\author{Benoit BERTRAND} 
\address{Departamento de \'Algebra\\
  Universidad Complutense de Madrid\\
 Avenida Complutense~s/n\\
 28040 Madrid\\ Espa\~na} 
\email{bertrand@mat.ucm.es} 
 \thanks{The author is supported by
  the European research network IHP-RAAG contract HPRN-CT-2001-00271}


\begin{abstract}
  A real algebraic variety is maximal (with respect to the Smith-Thom
  inequality) if the sum of the Betti numbers (with $\mathbb{Z}_2$
  coefficients) of the real part of the variety is equal to the sum of
  Betti numbers of its complex part.  We prove that there exist
  polytopes that are not Newton polytopes of any maximal hypersurface
  in the corresponding toric variety. On the other hand we show that
  for any polytope $\Delta$ there are families of hypersurfaces with
  the Newton polytopes $(\lambda\Delta )_{\lambda \in \mathbb{N}}$
  that are asymptotically maximal when $\lambda$ tends to infinity.
  We also show that these results generalize to complete
  intersections.

\end{abstract}

\maketitle

\newcommand\CQFD{\hfill $\Box$ \newline}
\newcommand\NN{{\mathbb{N}}}
\newcommand\ZZ{{\mathbb{Z}}}
\newcommand\PP{{\mathbb{P}}}
\newcommand\RR{{\mathbb{R}}}
\newcommand\CC{{\mathbb{C}}}
\newcommand\QQ{{\mathbb{Q}}}
\newcommand\bin[2]{{\mathrm{C}^{#2}_{#1}}}
\newcommand\FF{\mathcal{F}}
\newcommand\im{\operatorname{im }}
\newcommand\Vol{\operatorname{Vol}}
\newcommand\vol{\operatorname{vol}}
\newcommand\ra{$\RR A$ }
\newcommand\ca{$\CC A$ }
\newcommand\fl[2]{{\lfloor \frac{#1}{#2}\rfloor}}
\newcommand\cal[1]{\mathcal{#1}}


\newcommand{\N}{{\mathbb N}}
\newcommand{\A}{{\mathbb A}}
\newcommand{\Z}{{\mathbb Z}}
\newcommand{\R}{{\mathbb R}}
\newcommand{\C}{{\mathbb C}}
\newcommand{\Q}{{\mathbb Q}}
\newcommand{\Rpd}{{\mathbb R}P^2}
\newcommand{\Rpt}{{\mathbb R}P^3}
\newcommand{\Rpn}{{\mathbb R}P^n}
\newcommand{\Cpd}{{\mathbb C}P^2}
\newcommand{\Cpt}{{\mathbb C}P^3}
\newcommand{\Cpn}{{\mathbb C}P^n}
\newcommand{\ostatok}{A}
\newcommand{\mnote}{\marginpar}
\newcommand{\sign}{{\operatorname{sign}}}
\newcommand{\proofend}{\hfill$\Box$\bigskip}
\newcommand{\Tor}{{\operatorname{Tor}}}

\newcommand\newtext[1]{{\bf #1}}


\section{Introduction}

In 1876 Harnack showed that the maximal number of connected components
of a real algebraic plane projective curve of degree $m$ is
$(m-1)(m-2)/2 + 1$. He also proved that for any positive integer $m$
there exist curves of degree $m$ which are maximal in this sense (i.e.
with $(m-1)(m-2)/2 + 1$ connected components). Harnack's bound is
generalized to the case of any real algebraic variety by the
Smith-Thom inequality. Let $b_i(V;K)$ be the $i^{\mbox{\scriptsize
    th}}$ Betti number of a topological space $V$ with coefficients in
a field $K$ (i.e. $b_i(V;K)=\dim_K(H_i(V;K))$).  Denote by $b_*(V;K)$
the sum of the Betti numbers of $V$.  Let $X$ be a complex algebraic
variety equipped with an anti-holomorphic involution $c$. The real
part $\RR X$ of $X$ is the fixed point set of $c$.  Then the
Smith-Thom inequality states that $b_*(\RR X;\ZZ_2) \le b_*(X;\ZZ_2)$.
A variety $X$ for which $b_*(\RR X;\ZZ_2)= b_*(X;\ZZ_2)$ is called a
{\it maximal} variety or {\it $M$-variety}.  The question ``does a
given family of real algebraic varieties contain maximal elements?''
is one of the problems in topology of real algebraic varieties. For
the family of the hypersurfaces of a given degree in $\RR P^d$ a
positive answer is obtained in \cite{IteVir} using the combinatorial
Viro method called $T$-construction (see \cite{vir3}, \cite{vir6},
\cite{ite1}, and Theorem~\ref{IV}).
This question is, in general, a difficult problem.  Indeed we show
that Itenberg-Viro's theorem of existence of $M$-hypersurfaces of any
degree in the projective spaces of any dimension cannot be generalized
straightforwardly to all projective toric varieties.  More precisely,
in any dimension greater than or equal to $3$ there are polytopes
$\Delta$ such that no hypersurface in the toric variety $X_\Delta$
associated with $\Delta$, with the Newton polytope $\Delta$, is
maximal. However, in the $2$-dimensional case such a generalization of
the Harnack theorem holds (see Section~\ref{Mhypconterex}).

Let us first consider the $3$-dimensional case.  Let $k$ be a positive
integer number, and $\Delta_k$ be the tetrahedron in $\RR^3$ with
vertices $(0,0,0), (1,0,0), (0,1,0)$, and $(1,1,k)$.  Note that the
only integer points of $\Delta_k$ are its vertices.

\begin{figure}[tb]
\begin{center}
\input{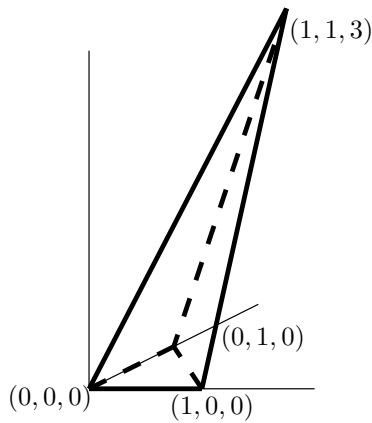}
\caption{Tetrahedron $\Delta_3$.}
\label{figure:fig1tetra}
\end{center}
\end{figure}

\begin{proposition}\label{skewtetra}
  For any odd $k \geq 3$ and any even $k \geq 8$, there is no
  maximal surface in $X_{\Delta_k}$ with the Newton polytope
  $\Delta_k$.
\end{proposition}

It is easy to generalize the above examples in dimension $3$ to higher
dimensions. From now on by {\it polytope} we mean a convex polytope
with integer vertices in the positive orthant $(\RR^+)^d = \{(x_1,
\ldots x_d) \in \RR^d \mid x_1 \geq 0, \ldots , x_d \geq 0\}$.

\begin{proposition}\label{skewsimplex}
  For any integer $d \ge 3$ there exist $d$-dimensional polytopes
  $\Delta$ such that no hypersurface in $X_\Delta$ with the Newton
  polytope $\Delta$ is maximal.
\end{proposition}

It is then natural to tackle the following weaker question.  Let
$\Delta$ be a $d$-dimensional polytope and $\{\lambda \cdot
\Delta\}_{\lambda\in\NN}$ the family of the multiples of $\Delta$.
Suppose that there exists a collection of polynomials
$\{P_\lambda\}_{\lambda\in\NN}$ satisfying the following conditions~:
\begin{enumerate}
\item the polytope $\lambda \cdot \Delta$ is the Newton polytope of
  $P_\lambda$,
\item
the total Betti numbers $b_*(\RR Z_\lambda;\ZZ_2)$ and
$b_*( Z_\lambda; \ZZ_2)$ are equivalent when $\lambda$ tends to
infinity
(here $Z_\lambda$ denotes the hypersurface in $X_\Delta$ defined by
$P_\lambda$).
\end{enumerate}

In this case we say that the family $\{Z_\lambda\}_{\lambda\in\NN}$ is
{\it asymptotically maximal}.  Given a $d$-dimensional polytope
$\Delta$ in $(\RR^+)^d$, does there exist an asymptotically maximal
family of hypersurfaces in $X_\Delta$?  A
positive answer to this question is given here.

\begin{theorem}\label{assmax}
  For any polytope $\Delta$ there exists an asymptotically maximal
  family of hypersurfaces $\{Z_\lambda\}_{\lambda\in\NN}$ in
  $X_\Delta$ such that for any $\lambda$ the Newton polytope of
  $Z_\lambda$ is $\lambda \cdot \Delta$.
\end{theorem}

The above statements have generalizations to complete intersections in
projective toric varieties. As a counterpart for Propositions \ref{skewtetra} and
 \ref{skewsimplex}  we show that, for any integer $d$ greater
than $2$ there exist polytopes $\Delta_d \subset (\RR^+)^d$ of
dimension $d$ such that the hypersurfaces defining a maximal complete
intersection in $X_{\Delta_d}$ cannot all have the Newton polytope
$\Delta_d$.

\begin{proposition}\label{cigencontrex}
  For any positive integers $d$ and $k$ such that $k\le d$ there exists a
  $d$-dimensional polytope $\Delta_d$ such that $k$ hypersurfaces
  defining a maximal complete intersection in $X_{\Delta_d}$ cannot
  all have the Newton polytope $\Delta_d$.
\end{proposition}

On the other hand the following theorem is a counterpart of
Theorem~\ref{assmax} for complete intersections.  Let $\Delta$ be a
$d$-dimensional polytope in $\RR^d$, and $k$ be an integer such that
$1 \le k \le d$. Knudsen-Mumford's theorem (see \cite{KKMS-D} p.161
and Theorem~\ref{KM}) asserts that there exists a positive integer $l$
such that $l \cdot \Delta$ admits a convex primitive triangulation.
Let $\lambda_1, \cdots ,\lambda_k$ be $k$ positive integers.  Denote
by $\Delta_{\lambda_i}$ the polytope $\lambda_i l \cdot \Delta$.  Let
$\{(\lambda_{1,m}, \cdots , \lambda_{k,m})\}_{m\in\NN}$ be a sequence
of $k$-tuples of positive integers such that $\lambda_{i,m}$ tends to
infinity for any $i = 1, \ldots , k$.  Let $\{(Z_{\lambda_{1,m}},
\cdots , Z_{\lambda_{k,m}})\}_{m}$ be a sequence of $k$-tuples of
algebraic hypersurfaces in $X_\Delta$ such that $Z_{\lambda_{i,m}}$
has the Newton polytope $\Delta_{\lambda_{i,m}}$. Assume that for any
natural number $m$ the variety $Y_m = Z_{1,m}\cap \cdots \cap Z_{k,m}$
is a complete intersection.

\begin{definition}
  Under the above hypotheses, the family $\{Y_m\}_{m\in\NN}$ is called
  {\it asymptotically maximal} if $b_*(\RR Y_m;\ZZ_2)$ is equivalent
  to $b_*(Y_m;\ZZ_2)$ when $m$ tends to infinity.
\end{definition}

\begin{theorem}\label{compassmaxth}
  Let $\Delta$ be a $d$-dimensional polytope, and $k$ be an integer
  number satisfying $1 \leq k \leq d$.  Let $\{(\lambda_{1,m}, \cdots
  , \lambda_{k,m})\}_{m\in\NN}$ be a sequence of $k$-tuples of natural
  numbers such that $\lambda_{i,m}$ tends to infinity for any $i = 1,
  \ldots , k$.  Then, there exists a sequence of $k$-tuples
  $\{(Z_{\lambda_{1,m}}, \cdots, Z_{\lambda_{k,m}})\}_{m\in\NN}$ of
  algebraic hypersurfaces in $X_\Delta$ such that

\begin{enumerate}
\item $Z_{\lambda_{i,m}}$ has the Newton polytope $\Delta_{\lambda_{i,m}}$
\item for any natural number $m$, the variety $Y_m = Z_{1,m}\cap \cdots \cap
  Z_{k,m}$ is a complete intersection,
\item the family $\{Y_m\}_{m\in\NN}$ is asymptotically maximal.
\end{enumerate}
\end{theorem}

{\bf Organization of the material.}  We first describe
  combinatorial patchworking and recall some results we will use.  In
  Section~\ref{hypassmax} we describe Itenberg and Viro construction
  of asymptotically maximal hypersurfaces in projective spaces. We
  then prove the existence of asymptotically maximal families of
  hypersurfaces for any Newton polytope (Theorem~\ref{assmax}).
  Proposition~\ref{skewtetra} and Proposition~\ref{cigencontrex} are
  proved respectively in Section~\ref{Mhypconterex} and in
  Section~\ref{Mcompintex}.  Finally, Section~\ref{compassmax} is
  devoted to the existence of asymptotically maximal families of
  complete intersections.  We describe there Itenberg and Viro
  construction of asymptotically maximal complete intersections in
  projective spaces and we prove Theorem~\ref{compassmaxth}.

  The author is grateful to Ilia Itenberg for his 
valuable advice.

\section{Preliminaries}\label{prelim}

\subsection{Toric varieties}

We fix here some conventions and notations, the construction of toric
varieties we use is based on the one described in \cite{ful}. Let
$\Delta$ be a polytope, $p$ a vertex of $\Delta$, and $\Gamma_1,
\cdots, \Gamma_k$ the facets of $\Delta$ containing $p$. To $p$ we
associate the cone $\sigma_p$ generated by the minimal integer inner
normal vectors of $\Gamma_1, \cdots, \Gamma_k$. The inner normal fan
$\mathfrak{E}_\Delta$ is the fan whose $d$-dimensional cones are the
cones $\sigma_p$ for all vertices $p$ of $\Delta$. The toric variety
$X_\Delta$ associated to $\Delta$ is the toric variety
$X(\mathfrak{E}_\Delta)$ associated to the fan $\mathfrak{E}_\Delta$
(see \cite{ful}).

\subsection{Combinatorial
patchworking}\label{tcon}

By a {\it subdivision} of a polytope we mean a
subdivision in convex polytopes (with integer vertices).
A subdivision
$\tau$ of  a  polytope $\Delta$ of dimension $d$
is called  {\it convex} if there exists a convex
piecewise-linear
function $\Phi : \Delta \to \RR$
whose domains of linearity coincide with the $d$-dimensional
polytopes
of $\tau$.

Let us  briefly describe the {\it combinatorial patchworking}, also called
{\it  $T$-construc\-tion}, which  is a particular
case  of  the Viro  method. A more detailed exposition can be found in
\cite{IteVir} (see also \cite{vir6} or \cite{GKZ}~p.~385).

Given a triple $(\Delta,\tau,D)$, where
$\Delta$ is a polytope, $\tau$ a convex triangulation of $\Delta$,
and $D$ a
distribution of signs at the vertices of $\tau$,
the combinatorial
patchworking, produces an algebraic hypersurface $Z$
in $X_\Delta$.

Let $\Delta$ be a $d$-dimensional polytope $(\RR^+)^d$ and $\tau$ be a
convex triangulation of $\Delta$.  Denote by $ s_{(i)}$ the reflection
with respect to the coordinate hyperplane $x_i=0 $ in $\RR^d$.
Consider the union $\Delta^*$ of all copies of $\Delta$ under the
compositions of reflections $ s_{(i)} $ and extend $\tau$ to a
triangulation $\tau^*$ of $\Delta^*$ by means of these reflections.
Let $D(\tau)$ be a sign distribution at the vertices of the
triangulation $\tau$ (i.e. each vertex is labelled with $+$ or $-$).
We extend $D(\tau)$ to a distribution of signs at the vertices of
$\tau^*$ using the following rule : for a vertex $a$ of $\tau^*$, one
has $\sign(s_{(i)}(a))=\sign(a)$ if the $i$-th coordinate of $a$ is
even, and $\sign(s_{(i)}(a))=-\sign(a)$, otherwise.

Let $\sigma$ be
a $d$-dimensional simplex of  $\tau^*$ with vertices of different
signs, and
$E$ be the hyperplane
piece which is the convex hull of the middle points of the
edges of $\sigma$
with endpoints of opposite signs. We separate vertices of
$\sigma$ labelled with $+$ from vertices labelled with $-$ by $E$.
The union of all these hyperplane pieces forms a piecewise-linear
hypersurface $H$.

For any facet $ \Gamma$ of $ \Delta^*$, let $N^\Gamma$ be a vector
normal to $\Gamma$. Let $F$ be a face of $\Delta^*$ and $\Gamma_1,
\dots, \Gamma_k$ be the facets containing $F$.  Let $L$ be the linear
space spanned by $N^{\Gamma_1}, \dots, N^{\Gamma_k}$. For any $v = (v_1,
\dots, v_d) \in L \cap \ZZ^d$ identify $F$ with 
${s_{(1)}}^{v_1} \circ
{s_{(2)}}^{v_2} \circ \dots \circ {s_{(d)}}^{v_d}(F)$.  Denote by $
\widetilde{\Delta}$ the result of the identifications.  The variety
$\widetilde{\Delta}$ is homeomorphic to the real part $\RR X_\Delta$
of $X_\Delta$ (see, for example, \cite{GKZ} Theorem~5.4 p. 383 or
\cite{stu} Proposition 2).

Denote by $\widetilde{H}$
the image of $H$ in $ \widetilde{\Delta}$.
Let $Q$ be a polynomial with the Newton polytope $\Delta$. It defines
a hypersurface $Z_0$ in the torus $(\CC^*)^d$ contained in
$X_\Delta$. The closure $Z$ of $Z_0$ in $X_\Delta$ is the hypersurface
defined by $Q$ in
$X_\Delta$. We call $\Delta$
the { \it Newton
polytope } of $Z$.

\begin{theorem}[T-construction, O. Viro]\label{Tcons}
Under the  hypotheses made above, there exists a
hypersurface $Z$ in $X_\Delta$ with
the Newton polytope $\Delta$
and a homeomorphism $h: \RR X_\Delta\to
\widetilde{\Delta}$ such that  $h(\RR  Z) = \widetilde{H}$.
\end{theorem}

The hypersurface  $Z$ in the above theorem
is called a  {\it real
algebraic $T$-hypersurface}.
A $d$-dimensional simplex with integer vertices
is called {\it primitive}
if its volume is equal to $\frac{1}{d !}$.
A triangulation $\tau$ of a $d$-dimensional polytope
is {\it primitive} if every $d$-simplex of the triangulation is primitive.
Let $\Delta$ be a $d$-dimensional  polytope. We call {\it lattice volume}
of $\Delta$ and  denote by $\Vol(\Delta)$
the volume normalized
so that a
primitive $d$-simplex has volume $1$.  The usual volume is denoted
by $\vol(\Delta)$.
If $\Delta$ is a $d$-dimensional polytope, then
$\Vol(\Delta) = d! \vol(\Delta)$.

\subsection{Sturmfels'
theorem
for complete intersections}

In \cite{stu} B. Sturmfels
proposed a combinatorial construction producing
complete intersections.
In fact, Sturmfels' construction is
an extended version of
the combinatorial patchworking.
We  quote here  this theorem in  the particular  case we
need. For the general statement and the proof we refer to \cite{stu}.

Let $\Delta_0$ be a $d$-dimensional polytope and $\lambda_1 ,\cdots
,\lambda_k$ positive integers, where $k \leq d$.  Denote by $\Delta_i$
the polytope $\lambda_i \cdot \Delta_0$ and by $\Delta$ the Minkowski
sum $\Delta_1 + \cdots +\Delta_k$.  Let $\nu_i$ be a piecewise-linear
convex function on $\Delta_i$ defining a triangulation $\tau_i$
with integer vertices.
For each
$\Delta_i$, choose a distribution of
signs $D_i$ at the  vertices of
$\tau_i$.

The initial data of the procedure of construction of a complete
intersection using Sturmfels' theorem are the polytopes $\Delta_i$,
the functions $\nu_i$ and the sign distributions $D_i$.  Apply the
$T$-construction for each triple $(\Delta_i, \tau_i, D_i)$ to
construct the hypersurfaces $S_{i}$. Let $D_i^*$ be the sign
distribution at the vertices of $\tau_i^*$.

  The functions $\nu_1 ,\cdots ,\nu_k$ define a convex
  decomposition of $\Delta$ in the following way (see \cite{stu},
  \cite{stu2} or \cite{Bih02}). Let $\bar{\Delta}_i$ be the convex
  hull of the set $\{(x,\nu_i(x)), x \in \Delta_i\}$ in $\RR^d \times
  \RR$ . Let $\bar{\Delta} \subset \RR^d \times \RR$ be the Minkowski
  sum $\bar{\Delta}_1 + \cdots + \bar{\Delta}_k$ and denote by $G$ the
  lower part of the boundary of $\bar{\Delta}$. Let $\nu$ be the
  piecewise-linear convex function of graph $G$ defined on $\Delta$
  (i.e. $G$ is the union of facets of $\bar{\Delta}$ whose inner
  normal vectors have positive last coordinate).  The function $\nu$
  defines a convex subdivision $\delta$ of $\Delta$ whose
  $d$-dimensional polytopes are the domains of linearity of $\nu$.
  Let $\Gamma$ be a polytope in $\delta$ and $\bar{\Gamma}$ its image
  by $\nu$. Then $\bar{\Gamma}$ can be uniquely written as the
  Minkowski sum $\bar{\Gamma}_1+ \cdots +\bar{\Gamma}_k$ where
  $\bar{\Gamma}_i$ is a face of $\bar{\Delta}_i$ for $i = 1,\cdots
  ,k$. This induces a decomposition of $\Gamma$ as a Minkowski sum
  $\Gamma = \Gamma_1 + \dots + \Gamma_k$ such that
  $\nu_i(\Gamma_i)=\bar{\Gamma}_i$.  
  Sturmfels' theorem requires the following genericity condition
  on the functions $\nu_i$.

\begin{definition}
  The $k$-tuple $\nu_1 ,\cdots ,\nu_k$ is said {\it sufficiently generic} if
  for any polytope $\Gamma$ of $\delta$, $\dim \bar{\Gamma} = \dim
  \bar{\Gamma}_1 + \dots + \dim \bar{\Gamma}_k$,  where $ \bar{\Gamma}
  = \bar{\Gamma}_1 + \dots + \bar{\Gamma}_k$ is the unique way to
  write $ \bar{\Gamma}$ as the Minkowski sum of faces of
  $\bar{\Delta}_1, \cdots, \bar{\Delta}_n$.
\end{definition}

We call { \it mixed subdivision} a subdivision $\delta$ obtained as
above from triangulations $\tau_1, \cdots, \tau_k$ and sufficiently
generic convex functions $\nu_1, \cdots, \nu_k$.  A mixed subdivision
$\delta$ is equipped with a decomposition of 
each of its polytopes
$\Gamma$ as a Minkowski sum $\Gamma = \Gamma_1 + \dots + \Gamma_k$,
where $\Gamma_i$ is a simplex of $\tau_i$. Two mixed subdivisions are
considered as equal if and only if they coincide as polyhedral
subdivisions, and each polytope of these subdivisions has the same
decomposition into a Minkowski sum in both of them.

Extend $\delta$ to a subdivision $\delta^*$ of $\Delta^*$ by means of
the reflections with respect to coordinate hyperplanes.  The extension
of the sign distribution to $\delta^*$ is as follows. Let $v$ be a
vertex of $\delta^*$, and let $v_1, \cdots ,v_k $ be the vertices of
$\tau_1^*, \cdots , \tau_k^*$ corresponding to $v$.  Then
$$\epsilon_j(v)=\sign(v_j).$$

For $j \in \{1, \cdots, k \}$ construct the hypersurface
$\widetilde{S_j}$ in the following way. For any polytope
$\Gamma^\prime$ in $\delta^*$, consider its symmetric copy $\Gamma$ in
$\delta$. There is a unique way to write $\Gamma = \Gamma_1 + \cdots
+\Gamma_k$ with $\Gamma_i$ in $\tau_i$ such that $\nu_1(\Gamma_1)+
\cdots+ \nu_k(\Gamma_k)= \nu(\Gamma)$. For $i \in \{ 1,\cdots, k \}$
let $\Gamma_i^\prime$ be the symmetric copy of $\Gamma_i$ in
$\tau_i^*$ such that $\Gamma^\prime = \Gamma_1^\prime + \cdots +
\Gamma_k^\prime$.  Define the hypersurface $S_j^*$ in $\Delta^*$ by
$S_j^*\cap \Gamma^\prime =\Gamma_1^\prime + \cdots + S_j \cap
\Gamma_j^\prime +\cdots + \Gamma_k^\prime $ for all $\Gamma^\prime$ in
$\delta^*$.  Let $\widetilde{S_j}$ be the image of $ S_j^*$ in
$\widetilde{\Delta}$.

\begin{theorem}[B. Sturmfels]\label{St}
With the above notation,  there exist hypersurfaces $Z_i$ with the Newton
polytopes $\Delta_i$, respectively,
and  a
homeomorphism
$f: \RR X_{\Delta} \to \widetilde{\Delta}$
such that
the hypersurfaces $Z_i$
define a complete intersection $Y$ in $X_\Delta$,
and $f$ sends $\RR Z_i$ (resp., $\RR Y$)
onto $\widetilde{S_i}$.
(resp., $\cap_{j=1 \cdots k}\widetilde{S_j} $).
\end{theorem}

\subsubsection{Cayley trick}\label{Cayley}
Instead  of constructing  the  complete intersection
in  the Min\-kowski sum of Newton polytopes, it
is convenient to use so-called Cayley trick (see, for example, \cite{stu2}).

Let $\Delta_1$, $\ldots$, $\Delta_k$
be convex polytopes with integer vertices
in $\R^d$ ($k \leq d$). For any $i = 1, \ldots , k$ put
$$\displaylines{\hat \Delta_i = \{(x_1, \ldots , x_{k+d}) \in \R^{k+d} \mid
x_i = 1; x_j = 0 \; \text{if} \; j \leq k \; \text{and} \; j \ne i; \cr
(x_{k+1}, \ldots , x_{k+d}) \in \Delta_i\}.}$$
The convex hull of $\hat \Delta_1$, $\ldots$, $\hat \Delta_k$ in $\R^{k+d}$
is called {\it Cayley polytope} and is denoted by $C(\Delta_1, \ldots , \Delta_k)$.
The intersection of $C(\Delta_1, \ldots , \Delta_k)$
with the subspace~$B \subset \R^{k+d}$
defined by $x_1 = \ldots = x_k = 1/k$ is naturally identified
with  the Minkowski sum~$\Delta$  of $\Delta_1$,  $\ldots$, $\Delta_k$
multiplied by $1/k$.
Thus, any triangulation of the Cayley polytope $C(\Delta_1, \ldots , \Delta_k)$
induces a subdivision of the Minkowski sum of $\Delta_1$, $\ldots$, $\Delta_k$.

The following lemma can be found, for example, in \cite{stu2}.

\begin{lemma}\label{Cayley-trick}
The correspondence described above establishes a bijection
between the set of convex triangulations with integer vertices of
$C(\Delta_1, \ldots , \Delta_k)$ and the set of mixed subdivisions
of the Minkowski sum of $\Delta_1$, $\ldots$, $\Delta_k$.\proofend
\end{lemma}

Denote by $C^*$ the union of the symmetric copies
of $C(\Delta_1, \ldots , \Delta_k)$ under
the reflections $s_{(i)}$, $i = k + 1, \ldots , k + n$,
where $s_{(i)}$ is the reflection of $\R^{k+d}$ with respect to
the hyperplane $\{x_i = 0\}$, and compositions of these reflections.

Choose a convex triangulation~$\tau$ of $C(\Delta_1, \ldots , \Delta_k)$ having
integer vertices and a distribution of signs at the vertices
of~$\tau$.
Extend the triangulation $\tau$
to a symmetric triangulation $\tau^*$ of
$C^*$
and the distribution of signs at the vertices of~$\tau$
to a distribution at the vertices of the extended triangulation by
the same rule as in Subsection~\ref{tcon}:
passing from a vertex to its
mirror image with respect to a coordinate hyperplane we preserve
its sign if the
distance from the vertex to the hyperplane is even, and change the sign if
the distance is odd.

For any ($k+d-1$)-dimensional simplex~$\gamma$ of~$\tau^*$
and any $j = 1, \ldots , k$
denote by~$\gamma_j$ the maximal face of~$\gamma$ which
belongs to a symmetric copy of $\hat \Delta_j$.
Let $K_j(\gamma)$ be the convex hull of
the middle points of the edges of $\gamma_j$ having endpoints of
opposite signs, and let $H(\gamma)$ be the intersection
of the join $K_1(\gamma)* \ldots * K_k(\gamma)$
with~$B$.
Denote by~$H$ the union of
the intersections $H(\gamma)$,
where $\gamma$ runs over all the ($k+d-1$)-dimensional simplices
of~$\tau^*$, and denote by $\widetilde H$ the image of~$H$
in $\widetilde{(\frac{1}{k}\Delta)}$.

The following statement is an immediate corollary
of Theorem~\ref{St}.

\begin{proposition}\label{proposition}
Assume that all the polytopes $\Delta_1$, $\ldots$, $\Delta_k$
are multiples of the same polytope~$\Pi$ with integer vertices.
Then, there exist nonsingular real hypersurfaces $Z_1$, $\ldots$, $Z_k$
in $X_\Pi$ with the Newton polytopes $\Delta_1$, $\ldots$, $\Delta_k$, respectively,
and a homeomorphism $f: \R X_\Pi \to \widetilde{(\frac{1}{k}\Delta)}$
such that the hypersurfaces $Z_1$, $\ldots$, $Z_k$ define a complete
intersection~$Y$ in $X_\Pi$ and~$f$ maps
the set of real points $\R Y$ of $Y$ onto $\widetilde H$.\proofend
\end{proposition}

\subsection{Formulae for the Betti numbers}

V. Danilov and A. Khovanskii \cite{DanKho} computed the
Hodge numbers of a smooth hypersurface in a toric variety $X_\Delta$ in
terms of the polytope $\Delta$ involving in particular the
coefficients of the Ehrhart polynomial of $\Delta$
 (see \cite{Ehr} or \cite{Ehr67}).
Our aim being to investigate asymptotical behaviors of certain families
of hypersurfaces or complete intersections, we need only the simpler results
that are quoted below.

\begin{definition}
A $d$-dimensional polytope $\Delta$  is simple if for each vertex
$a$ of $\Delta$, the number of edges of $\Delta$ containing $a$ is $d$.
\end{definition}

Let $l^*(\Delta)$ be the number of integer points in the interior of
$\Delta$ (i.e. $l^*(\Delta)= \# (\ZZ^d \cap (\Delta \setminus \partial
\Delta))$ ). The following statement can be found in \cite{DanKho}
Section~{5.11}.

\begin{lemma}\label{bZ}
Let $\Delta$ be a $3$-dimensional simple polytope,
and $Z$ be an algebraic hypersurface of
$X_\Delta$  with the Newton polytope  $\Delta$.  Then $b_*(Z;\CC)  =l^*(2
\Delta)   -2  l^*(\Delta)   -  \sum_{\Gamma   \in  \FF_2   (\Delta)  }
(l^*(\Gamma)-1)-1$.
\end{lemma}

The  following  two  propositions  can be  derived  from  Khovanskii's
results (see \cite{Kho1} and \cite{Kho2}) or can be found in \cite{Mik2}.

\begin{proposition}\label{DKhypass}
Let $\Delta$ be a polytope,
and $\{Z_\lambda\}_{\lambda\in\NN}$ be a family of
algebraic hypersurfaces in $X_\Delta$
with the Newton polytopes $\lambda \cdot \Delta$.
Then
$b_*(Z_\lambda;\ZZ_2)$ is equivalent to $\Vol(\lambda \cdot \Delta)$
when $\lambda$ tends to
infinity.
\end{proposition}

Denote by $\Vol(\Delta_1  , \cdots ,  \Delta_k)$
the mixed volume of
the  polytopes   $\Delta_1  ,  \cdots  ,  \Delta_k$.
We choose a normalization of  the mixed volume
in such a way that for a  primitive simplex
$\sigma$ we have $\Vol(\sigma, \cdots ,\sigma)= 1$.

\begin{proposition}\label{DKcompass}
Let $\Delta$ be a
$d$-dimensional polytope, and $k$ be a positive integer
satisfying $k \leq d$.
Assume that for any collection
$\lambda_1   ,  \cdots   ,  \lambda_k$ of positive integers
we have a collection of $k$ hypersurfaces $Z_{\lambda_{1}},   \cdots,
Z_{\lambda_{k}}$ in $X_\Delta$ with the Newton polytopes
$\lambda_1 \cdot \Delta, \ldots ,\lambda_k \cdot \Delta$,
respectively, such that $Z_{\lambda_{1}},   \cdots,    Z_{\lambda_{k}}$
define a complete intersection $Y_{\lambda_1 , \ldots , \lambda_k}$
in $X_\Delta$.
Then
$b_*(Y_{\lambda_1, \ldots, \lambda_k}; \ZZ_2)$
is equivalent
to $\Vol(\lambda_1 \cdot \Delta,  \cdots, \lambda_k \cdot \Delta)$
when $\lambda_i$
tends to infinity for all $i$.
\end{proposition}

We  also use the  following result of  Khovanskii on  the Euler
characteristic of a complete intersection
in the torus $(\CC^*)^d $ (see \cite{Kho2}).

\begin{theorem}[A. Khovanskii]\label{Kho}
Let $Y$ be a complete intersection in $(\CC^*)^d $ defined by
polynomials $P_1 ,\cdots ,P_k$ with the Newton polytopes $\Delta_1,
\cdots ,\Delta_k$, respectively. Then, the Euler characteristic of
$Y$ is the homogeneous  term of  degree $d$  of
$$\Delta_1(1  +
\Delta_1)^{-1}\cdot \> \cdots \> \cdot \Delta_k(1
+\Delta_k)^{-1},$$
where the  product of  $d$ polytopes stands for
their mixed volume and $(1 + \Delta_i)^{-1}$ stands for the series
$\sum_{j=0}^\infty (-1)^j (\Delta_i)^j$.
\end{theorem}

In the case of two $3$-dimensional polytopes we use the following
direct consequence of
Theorem~\ref{Kho}.

\begin{bbcor}\label{CorKho}
Let  $\Delta$  be  a  simple  $3$-dimensional  polytope
and  $\lambda_1$  and
$\lambda_2$ be positive integers. For $i=1,2$ put
$\Delta_i=\lambda_i \cdot \Delta$.
Let $Y$ be a complete  intersection in
$X_\Delta$ defined by polynomials $P_1$ and $P_2$ with the Newton polytopes
$\Delta_1$  and  $\Delta_2$,  respectively. Then,
$b_*(Y;\CC)=(\lambda_1^2
\lambda_2  + \lambda_2^2  \lambda_1)
\Vol(\Delta)
-  \sum_{\Gamma \in
\mathcal{F}_2(\Delta)}
\lambda_1\lambda_2
\Vol(\Gamma)
+ 4$.
\end{bbcor}

\proof By Theorem~\ref{Kho}, the Euler characteristic $\chi(Y)$ of $Y$
is given by $\chi(Y)=-(\lambda_1^2\lambda_2 + \lambda_2^2 \lambda_1)
\Vol(\Delta) + \sum_{\Gamma \in \mathcal{F}_2(\Delta)}
\lambda_1\lambda_2 \Vol(\Gamma)$.  Since $b_*(Y;\CC) = -\chi(Y) + 4$,
we have the desired result. \CQFD

\section{Asymptotically
maximal families of hypersurfaces}\label{hypassmax}

\subsection{Auxiliary statements}\label{const}

This section is devoted to the proof of Theorem~\ref{assmax} on existence of
asymptotically maximal families of hypersurfaces.
The proof
is based on two important results.

In \cite{IteVir} I. Itenberg and O. Viro, using the
$T$-construction, proved that there
exist M-hypersurfaces
of any degree in the projective
space of any dimension.

\begin{theorem}[I. Itenberg and O. Viro]\label{IV}
Let $d$ and $m$ be natural numbers,
and $T_1^d$ be a primitive $d$-dimensional simplex.
Put $T_m^d= m \cdot T_1^d$.
Then, there exists a
primitive convex
triangulation ${\tau}_{T_m^d}$  of   $T_m^d$ and a sign
distribution $D({\tau}_{T_m^d})$
at the vertices of ${\tau}_{T_m^d}$
such that the $T$-hypersurface $Z_d^m$ obtained via
the combinatorial patchworking from
${\tau}_{T_m^d}$ and $D({\tau}_{T_m^d})$ is maximal.
\end{theorem}

The second important result we
use is due to  F. Knudsen and D. Mumford \cite{KKMS-D}.

\begin{theorem}[F. Knudsen and D. Mumford]\label{KM}
Let $\Delta$ be a polytope. There exists a positive integer $l$ such
that $l \cdot \Delta$ admits a convex primitive triangulation.
\end{theorem}

In the sequel, when there is no ambiguity on the triangulation of a
polytope $\Delta$ and the sign distribution chosen, we denote by
$H_\Delta$ the piecewise-linear hypersurface in $\Delta_*$ obtained by
$T$-construction, $\widetilde{H_\Delta}$ its image in
$\widetilde{\Delta}$, and $Z_\Delta$ the corresponding hypersurface in
$X_\Delta$.

\subsection{Itenberg-Viro asymptotical construction}\label{IVhypass}

In fact, we use only the following asymptotical version of
Theorem~\ref{IV}.

\begin{theorem}[I. Itenberg and O. Viro]\label{IVtheorem1}
  For any positive integers $m$ and $d$ such that $m \geq d + 1$,
  there exists a hypersurface~$X$ of degree $m$ in $\RR P^d$ such that
$$b_*(\R X;\ZZ_2) \geq (m - 2)(m - 3) \ldots (m - d - 1).$$
\end{theorem}

The proof of this asymptotical version is much simpler than the proof
of Theorem~\ref{IV}.  It can be extracted from \cite{IteVir} and was
communicated to us by the authors of~\cite{IteVir}.  We reproduce
their proof here for the completeness.
\subsubsection{Proof of Theorem~\ref{IVtheorem1}}\label{IVsection1}

We describe a triangulation $\tau$ of the standard
simplex $T = T^d_m$
and a distribution of signs at the integer points of $T$
which provide via the combinatorial patchworking theorem
a hypersurface with the properties formulated
in Theorem~\ref{IVtheorem1}.

To construct the triangulation $\tau$, we use induction on $d$.  If $d
= 1$, the triangulation of $[0, m]$ is formed by $m$ intervals $[0,
1]$, \dots, $[m - 1, m]$ for any $m$.
Assume that for all natural $k<d$ the triangulations of the standard
$k$-dimensional simplices of all sizes are constructed and consider
the $d$-dimensional one of size $m$.

Denote by $x_1, \ldots , x_d$ the coordinates in $\R^d$.  Let $T_j^{d
- 1} = T \cap \{x_d = m - j\}$ and $T_j$ be the image of $T_j^{d - 1}$
under the orthogonal projection to the coordinate hyperplane $\{x_d =
0\}$.  Numerate the vertices of each simplex $T_1, \ldots
, T_{m - 1}, T_m= T_m^{d - 1}$ as follows: assign $1$ to the
vertex at the origin and $i+1$ to the vertex with nonzero
coordinate at the $i$-th place.  Assign to the vertices
of $T_1^{d - 1}, \ldots , T_{m - 1}^{d - 1}$ the numbers of
their projections.
A  triangulation of each
simplex $T_0, \ldots , T_{m - 1}$ is constructed.
Take the corresponding triangulations in the simplices
$T_j^{d - 1}$.

Let $l$ be a nonnegative integer not greater than $d - 1$.
If $m - j$ is even, denote by $T_j^{(l)}$ the $l$-face of $T_j^{d - 1}$
which is the convex hull of the vertices with numbers
$1, \ldots , l + 1$.
If $m - j$ is odd denote by $T_j^{(l)}$ the $l$-face of $T_j^{d - 1}$
which is the convex hull of the vertices with numbers
$d - l, \ldots , d$.

Now for any integer $0 \leq j \leq m - 1$ and any integer $0 \leq l
\leq d - 1$, take the join $T_{j + 1}^{(l)} * T_j^{(d - 1 - l)}$.  The
triangulations of $T_{j + 1}^{(l)}$ and $T_j^{(d - 1 - l)}$
define a triangulation of $T_{j + 1}^{(l)} * T_j^{(d - 1 - l)}$.
This gives rise to the desired triangulation $\tau$ of $T$.
One can see that $\tau$ is convex.

The distribution of signs at the vertices of $\tau$
is given by the following rule.
The vertex gets the sign ``$+$'' if the sum of its coordinates
is even, and it gets the sign ``$-$'' otherwise.

\begin{lemma}\label{IVlemma} For the hypersurface~$X$ of degree
$m$ in $\RR P^d$ provided according to the combinatorial patchworking
theorem by the
triangulation $\tau$ and the distribution of signs defined above,
one has
$$b_*(\R X;\ZZ_2) \geq \begin{cases}
(m - 2)(m - 3) \ldots (m - d - 1), &\text{if
$m \geq d + 1$,}\\
\; 0,
&\text{otherwise.}
\end{cases}$$
\end{lemma}

To prove Lemma~\ref{IVlemma}
we define a collection of cycles $c_i$, $i \in I$
of $\widetilde H$ (in fact, any $c_i$
is also a cycle of the hypersurface
$H \subset T_*$, and moreover, of the hypersurface
$H \cap (\RR^*)^d$).
The cycles $c_i$ are called {\em narrow}.

The collection of narrow cycles $c_i$ is constructed together with
a collection of
{\it dual cycles} $b_i$.
Any dual cycle $b_i$ is a $(d - 1 - p)$-cycle
in $\widetilde T \setminus \widetilde H$
(where $p$ is the dimension of $c_i$)
composed by simplices of $\tau_*$ and representing
a homological class such that its
linking number with any $p$-dimensional narrow
cycle $c_k$ is $\delta_{ik}$.

Let us fix some notations.
For any simplex $T_j^{(l)}$ (where $1 \leq j \leq m$
and $0 \leq l \leq d - 1$), denote by $(T_j^{(l)})_*$
the union of the symmetric copies of $T_j^{(l)}$ under the reflections
with respect to coordinate hyperplanes
$\{x_i = 0\}$, where $i = 1, \ldots, l$, if $m - j$ is even,
and $i = d - l, \ldots , d - 1$, if $m - j$ is odd, and compositions
of these reflections.

Any simplex $T_j^{(l)}$ is naturally identified with
the standard simplex $T^l_j$ in $\R^l$ with vertices
$(0, \ldots , 0)$, $(j, 0, \ldots , 0)$,
$\ldots, (0, \ldots , 0, j)$ via the linear
map ${\cal L}_j^l: T_j^{(l)} \to T^l_j$ sending
\begin{enumerate}
\item
the vertex with number $i$ of $T_j^{(l)}$ to the vertex of $T^l_j$
with the same number, if $m - j$ is even,
\item
the vertex with number $i$ of $T_j^{(l)}$
to the vertex of $T^l_j$
with the number $i - d + l + 1$, if $m - j$ is odd.
\end{enumerate}
It is easy to see that ${\cal L}_j^l$ is simplicial
with respect to the chosen triangulations of $T_j^{(l)}$ and $T^l_j$.
The natural extension of ${\cal L}_j^l$ to $(T_j^{(l)})_*$
identifies $(T_j^{(l)})_*$ with $(T^l_j)_*$
and respects the chosen triangulations.

By {\it a symmetry} we mean a composition of reflections with respect
to coordinate hyperplanes.
Let $s_{(i)}$ be the reflection of $\R^d$ with respect to
the hyperplane $\{x_i = 0\}$, $i = 1, \ldots , d$.
Denote by $s_j^l$ the symmetry of $(T_j^{l + 1})_*$ which is
identical if $m - j$ is even, and coincides with the restriction of
$s_{(d - l - 1)} \circ \ldots \circ s_{(d - 1)}$ on $(T_j^{l + 1})_*$
if $m - j$ is odd.

The narrow cycles and their dual cycles
are defined below using induction
on $d$. For $d = 1$ the narrow cycles are the pairs of points
$$(1/2, 3/2), \ldots , ((2m - 5)/2, (2m - 3)/2).$$
The dual cycles are pairs of vertices
$$(1, m - 1), (2, m), (3, m + 1), \ldots , (m - 2, m),$$
if~$m$ is even, and pairs of vertices
$$(1, m), (2, m - 1), (3, m), \ldots , (m - 2, m),$$
if~$m$ is odd.

Assume that for all natural $m$ and all natural $k<d$
the narrow cycles $c_i$ in the hypersurface
$\widetilde H \subset \widetilde T_m^k$
and the dual cycles $b_i$ in $\widetilde T_m^k \setminus \widetilde H$
are constructed.
The narrow cycles of the hypersurface in $\widetilde T_m^d$
are divided into $3$ families.

\vskip10pt

{\bf Horizontal Cycles. } The initial data for constructing a cycle
of the first family consist of an integer $j$ satisfying inequality
$1 \leq j \leq m - 1$ and a narrow cycle of the hypersurface in
$T^{d-1}_*$ constructed at the previous step.
In the copy $(T_j^{d - 1})_*$ of $T^{d-1}_*$,
take the copy $c$ of this cycle and $b$ of
its dual cycle.

There exists exactly one symmetric copy of $T_{j + 1}^0$ incident to
$b$.  It is $T_{j + 1}^0$ itself, if $m - j$ is odd, and either $T_{j
+ 1}^0$, or $s_{(d - 1)}(T_{j + 1}^0)$, if $m - j$ is even.  If the
sign of the symmetric copy $s(T_{j + 1}^0)$ of $T_{j + 1}^0$ incident
to $b$ is opposite to the sign of $c$, we include $c$ in the
collection of narrow cycles of $\widetilde H$.  Otherwise take
$s_{(d)}(c)$ as a narrow cycle of $\widetilde H$.  The dual cycle of $c$
(resp., $s_{(d)}(c)$)
is the suspension of~$b$ (resp., $s_{(d)}(b)$)
with the vertex $s(T_{j + 1}^0)$ (resp., $s_{(d)}(s(T_{j + 1}^0))$)
and with the vertex $s(T_{j - 1}^0)$
(resp., $s_{(d)}(s(T_{j - 1}^0))$).

\vskip10pt

{\bf Co-Horizontal Cycles.  } The initial data for constructing a
cycle of the second family are the same as in the case of the
horizontal cycles: the data consist of an integer $j$ satisfying
inequality $1 \leq j \leq m - 1$ and a narrow cycle of the
hypersurface in $T^{d-1}_*$.

In the copy $(T_j^{d - 1})_*$ of $T^{d-1}_*$, take the copy $c$ of
this cycle and $b$ of its dual cycle. If the sign of the symmetric
copy $s(T_{j + 1}^0)$ of $T_{j + 1}^0$ incident to $b$ coincides with
the sign of $c$, take $b$ as dual cycle of a narrow cycle of
$\widetilde H$.  Otherwise take $s_{(d)}(b)$.  The corresponding
narrow cycle is a suspension of~$c$ (resp., $s_{(d)}(c)$).

\vskip10pt

{\bf Join Cycles. } The initial data consist
of integers~$j$ and~$l$
satisfying inequalities $1 \leq j \leq m - 1$,
$1 \leq l \leq d - 2$,
the copy $c_1\subset (T_{j + 1}^l)_*$ of a narrow cycle
of the hypersurface in $(T^l_{j + 1})_*$, the copy
$c_2\subset (T_j^{d - 1 - l})_*$ of a narrow cycle
of the hypersurface in $(T^{d - 1 - l}_j)_*$ and the copies
$b_1\subset (T_{j + 1}^l)_*$ and $b_2\subset (T_j^{d - 1 - l})_*$ of
the dual cycles of these narrow cycles.

One of the joins $b_1 * b_2$ and $s_{j + 1}^l(b_1) * s_j^{d - 1 -
l}(b_2)$, belongs to $\tau_*$; denote it by $J$.  If the signs of
$c_1$ and $c_2$ coincide, take $J$ as the dual cycle of a cycle of
$\widetilde H$.  Otherwise take $s_{(d)}(J)$.  The corresponding
narrow cycle is either $c_1*c_2$, or
$s_{j + 1}^l(c_1) * s_j^{d - 1 - l}(c_2)$,
or $s_{(d)}(c_1 * c_2)$, or
$s_{(d)}(s_{j + 1}^l(c_1) * s_j^{d - 1 - l}(c_2))$.

\vskip10pt

{\bf Proof of Lemma~\ref{IVlemma}.}
Both $c_i$ and $b_i$ with $i\in I$ are $\Z_2$-cycles
homologous
to zero in $\widetilde{T}$, which is homeomorphic to the projective
space of dimension $d$. The sum of dimensions of $c_i$ and
$b_i$ is $d-1$. Thus we can consider the linking
number of $c_i$ with $i\in I$ and $b_k$, $k\in I$
taking values in
$\Z_2$. Each $c_i$ bounds an obvious ball in $\widetilde T$. This ball
meets $b_i$ in a single point transversally and is disjoint
with $b_k$ for $k\ne i$ and $i,k\in I$. Hence the  linking
number of $c_i$ and $b_k$ is $\delta_{ik}$.

Therefore  the collections of homology classes realized in
$\widetilde T \setminus \widetilde H$ and $\widetilde H$ by
$b_i, i \in I$ and  $c_i, i \in I$, respectively,
generate subspaces of
$H_*(\widetilde T \setminus \widetilde H;\Z_2)$
and $H_*(\widetilde H;\Z_2)$ and are
dual bases of the subspaces with respect to the restriction of the
Alexander duality.  Hence $c_i$ with $i\in I$
realize linearly independent
$Z_2$-homology classes of $\widetilde H$.

It remains to show that the number of narrow cycles is at least
$$(m - 2)(m - 3) \ldots (m - d - 1),$$
if $m \geq d + 1$.
The statement can be proved by induction on~$d$.
The base $d = 1$ is evident. To prove the induction step
notice, first, that the statement is evidently true for
$m = d + 1$. Now, we use the induction
on~$m$ and obtain the required statement
from the inequality
$$\displaylines{
(m - 3)(m - 4) \ldots (m - d - 2) + 2(m - 3)(m - 4)
\ldots (m - d - 1) \cr
+ \sum_{k = 1}^{d - 2}\left[(m - 2)(m - 3) \ldots (m - k - 1)\right]
\left[(m - 3)(m - 4) \ldots (m - d + k - 1)\right] \cr \geq
(m - 2)(m - 3) \ldots (m - d - 1).}$$
This finishes the proofs of
Lemma~\ref{IVlemma} and Theorem~\ref{IVtheorem1}.
\proofend

\begin{remark}
The family of hypersurfaces
in $\RR P^d$ constructed
in Theorem~\ref{IVtheorem1}
is asymptotically maximal.
\end{remark}

\proof
Indeed, the total Betti
number of a nonsingular hypersurface
of degree $m$ in $\CC P^d$
is equal to
$\frac{(m-1)^{d+1} - (-1)^{d+1}}{m} +  d + (-1)^{d+1}$.
This number
is equivalent to  $ (m - 2)(m - 3) \ldots (m -  d - 1)$ when $m$
tends
to
infinity. \CQFD

\subsection{Proof of Theorem~\ref{assmax}}\label{assmax-section}

For a positive integer $\lambda$ put
$\Delta_\lambda = \lambda \cdot \Delta$. Let $l$ be
a positive integer such
that $\Delta_l$ admits a primitive convex triangulation $\tau$
(see Theorem~\ref{KM}).
Denote by $\nu$
a function certifying the convexity of $\tau$.
Let $\tau_\lambda$ be the triangulation of $\Delta_{\lambda l}$
obtained from $\tau$ by
multiplication of its simplices by $\lambda$.

We can assume that $\lambda > d + 1$.
Let
$\delta$ be a
$d$-dimensional simplex of $\tau$.
The convex hull of the interior integer points
of $\lambda \cdot \delta$ is a
$d$-dimensional simplex $(\lambda-(d+1)) \cdot \delta$.
Put $\delta_\lambda=\lambda
\cdot \delta$ and $\delta_\lambda^\prime=(\lambda-(d+1)) \cdot \delta$.
For any
$d$-dimensional simplex $\delta_\lambda$ of
$\tau_\lambda$, apply the construction
of Lemma~\ref{IVlemma} to the
convex  hull $\delta_\lambda^\prime$ of
the interior integer points of $\delta_\lambda$.
Complete the triangulation of
$\delta_\lambda^\prime$
to a convex triangulation of $\delta_\lambda$ whose only
extra vertices are the vertices of $\delta_\lambda$ in
the following way. Let $\nu_{\lambda - (d+1)}$ be
a convex piecewise-linear function
certifying the convexity of the triangulation of
$\delta_\lambda^\prime$.
Define a convex function $\nu_\lambda^\delta$ on  $\delta_\lambda$
choosing the values of $\nu_{\lambda - (d+1)}$
at the integer points of $\delta_\lambda^\prime$ and
the
value $v$ at the vertices
of
$\delta_\lambda$, where $v$ is large enough
(the graph of $\nu_\lambda^\delta$ is the lower part of the convex hull of the defined points in $\delta_\lambda \times \RR$).
Note that $\nu_\lambda^\delta$ restricted to $\delta_\lambda^\prime$
coincides with $\nu_{\lambda  - (d+1)}$.
If the decomposition defined by
$\nu_\lambda^\delta$ is not a triangulation, we
slightly perturb $\nu_{\lambda - (d+1)}$
(without changing the triangulation of $\delta_\lambda^\prime$)
to
break the
polytopes of the subdivision which are not simplices.
Denote by $\tau_\lambda^\delta$ the obtained triangulation
of $\delta_\lambda$.

The only vertices of $\tau_\lambda^\delta$ in $\delta_\lambda
\setminus \delta_\lambda^\prime$ are the vertices of $\delta_\lambda$.
One can choose the same value $v$ of the functions $\nu_{\lambda}^\delta$
at the vertices of
all the
$d$-dimensional simplices $\delta$ of $\tau_\lambda$.
Hence, the functions
$\nu_\lambda^\delta$ can be glued together to
form a
piecewise-linear function $\nu_\lambda$
on $\Delta_{\lambda l}$ which is,
by construction, convex on each
$d$-dimensional simplex of $\tau_\lambda$.
Let $\nu^\prime$ be a function certifying the convexity of
$\tau_\lambda$.
Then, for sufficiently small $\epsilon > 0$ the function
$\nu=\nu^\prime  + \epsilon \nu_\lambda$
certifies the convexity of the triangulation obtained by gluing the
triangulations of the $d$-dimensional simplices of $\tau_\lambda$.
Thus, one gets a convex
triangulation $\tau_\lambda^l$ of $\Delta_{\lambda l}$.
Choose a sign  distribution $D(\tau_\lambda  ^l)$
at the vertices of $\tau_\lambda^l$
in such a way that on each
simplex $\delta_\lambda^\prime$
the distribution
coincides with the one
Lemma~\ref{IVlemma}.
Let $Z_{\Delta_{\lambda l}}$ be the hypersurface obtained
via the combinatorial patchworking
from $\tau_\lambda^l$ and $D(\tau_\lambda^l)$.

\begin{proposition}
The  family of hypersurfaces  $Z_{\Delta_{\lambda l}}$ of
$X_\Delta$ constructed above is asymptotically maximal.
\end{proposition}

\proof
The total Betti number of $Z_{\Delta_{\lambda l}}$ is equivalent
to $\Vol(\Delta_{\lambda l})$ when $\lambda$ tends to infinity
(see Proposition~\ref{DKhypass}).
For each $d$-dimensional simplex $\delta$ of $\tau_\lambda$
consider the narrow cycles of $H_{\Delta_{\lambda l}}\cap
(\delta_\lambda^\prime)_*$ which are constructed in the proof of
Lemma~\ref{IVlemma}. Since the narrow cycles are constructed with
the dual cycles, the union of the obtained collections of narrow
cycles consists of linearly independent cycles. Thus, $b_*(\RR
Z_{\Delta_{\lambda l}};\ZZ_2) \ge \Vol(\Delta_{l}) n_\lambda$, where
$n_\lambda$ is the  number of narrow cycles  in each
$\delta_\lambda^\prime$. Since $n_\lambda \sim
\Vol(\delta_\lambda^\prime)$, we have
$n_\lambda\sim\Vol(\delta_\lambda)$. So, $b_*(\RR
Z_{\Delta_{\lambda l}};\ZZ_2)$ is equivalent to $ \Vol(\Delta_{l})
\Vol(\delta_\lambda)$. The latter number is equal to
$\Vol(\Delta_{\lambda l})$. \CQFD

\section{Newton polytopes without maximal hypersurfaces}\label{Mhypconterex}

Before giving the proof Proposition~\ref{skewtetra} let us
consider the lower dimensional cases.
Clearly, if $\Delta$ is an interval $[a, b]$ in $\RR$,
where $a$ and $b$
are nonnegative integers, then there exists
a maximal $0$-dimensional subvariety in $\CC P^1 = X_\Delta$
with the Newton polygon $\Delta$.

If $\Delta$ is a polygon in the first quadrant of $\RR^2$,
then again there exists a maximal curve in $X_\Delta$
with the Newton polygon $\Delta$.
Such a curve can be constructed by the combinatorial patchworking:
it suffices to take as initial data
a primitive convex triangulation of $\Delta$
equipped with the following distribution of signs:
an integer point $(i, j)$ of $\Delta$ gets
the sign ``-'' if $i$ and $j$ are both even,
and gets the sign ``+'', otherwise
(see for example  \cite{Iterag}, \cite {IteVir2}, \cite{Haa}).

\vskip10pt

{\bf Proof of Proposition~\ref{skewtetra}.}
The proof of Proposition~\ref{skewtetra} relies on the
estimation of the Betti numbers of the
complex
and real parts
of a real algebraic  surface  $Z_k$ in $X_{\Delta_k}$
with the
Newton polytope
$\Delta_k$. The Betti numbers $b_*(Z_k; \CC)$ are given
by Lemma~\ref{bZ}.
We have $b_*  ( Z_k; \CC)= l^*(2  \Delta_k) -2 l^*(\Delta_k)
- \sum_{\Gamma \in \FF_2 (\Delta_k) }
(l^*(\Gamma)-1)-1$.
Since $l^*(2   \Delta_k)=  k-1$  and  $l^*(\Delta_k)=0$,
we get  $b_* (  Z_k; \CC)  = k+2$.
Thus, $b_* ( Z_k;\ZZ _2 ) \ge k+2$.

To estimate $b_* (\RR  Z_k;\ZZ _2 )$ we
consider
two cases.
If $k$ is odd, $\Delta_k$
is an elementary tetrahedron, and
$\RR Z_k$
is homeomorphic to the projective plane.
Thus, in this case,
$b_* (\RR Z_k;\ZZ _2)= 3$.

  If $k$ is even, $\Delta_k$ has either $6$ or $8$ nonempty
  symmetric copies.  In the first case $\RR Z_k$ is homeomorphic to
  three spheres with some points identified.  Each of the spheres has
  $4$ marked points. Pairs of marked points are identified in the
  following way.  Two marked points of each sphere are identified with
  two marked points of another sphere, and the two other marked points
  are identified with the marked points of the remaining sphere. Then
  the Euler characteristic is zero and $b_* (\RR Z_k;\ZZ _2) = 8$.  In
  the case of $8$ nonempty symmetric copies, $\RR Z_k$ is homeomorphic
  to four spheres with some points identified. Each sphere has three
  marked points. Pairs of marked points are identified in the
  following
  way: on each sphere the three marked points are identified with
  marked points of three different spheres. Thus the Euler
  characteristic is $2$ and we also have $b_* (\RR Z_k;\ZZ _2) = 8$.

  Thus, for $k$ even greater than or equal to $8$ and for $k$ odd
  greater than or equal to $3$, there is no maximal surface in
  $X_{\Delta_k}$ with the Newton polytope $\Delta_k$.

\vskip10pt

{\bf Proof of Proposition~\ref{skewsimplex}.}
Fix an integer $d \geq 3$ and consider
a family
$\{\sigma_k\}_{k \in \NN}$ of $d$-dimensional simplices
in $\RR^d$ such that
their vertices
are their only integer points
and $\Vol(\sigma_k)=k$. For example, one can take for $\sigma_k$
the simplex in $\RR^d$ with vertices
$$\displaylines{(0, 0, \ldots , 0, 0),\> (1, 0, \ldots , 0, 0),\>
(0, 1, \ldots , 0, 0),\>
\ldots ,
(0, 0, \ldots , 1, 0),\cr
\text{\rm and } (1, 1, \ldots , 1, k).}$$

Let $Z_k$ be any hypersurface in $X_{\sigma_k}$.
By Proposition~\ref{DKhypass} $b_*(Z_k;\CC)$ tends
to infinity when
$k$ does, and so does $b_*(Z_k;\ZZ_2)$. Meanwhile, $b_*(\RR Z_k;\ZZ_2)$
is bounded (for example,
by the number of simplices in $\sigma_k^*$).
So there exists a number $k_0$ such
that for any integer $k > k_0$ and any hypersurface $Z_k$ in $X_{\sigma_k}$
one has $b_*(\RR Z_k;\ZZ_2) < b_*(Z_k;\ZZ_2)$.\CQFD

\section{Newton polytopes without maximal complete intersection}\label{Mcompintex}

Let us first consider the case of complete intersections of two surfaces.
Let $\Delta_k$ be the tetrahedron in $\RR^3$ with vertices $(0,0,0),
(1,0,0), (0,1,0)$ and $(1,1,k)$.

\begin{proposition}\label{ciconterex}
  Let $k \geq 5$ be an integer, and $Z_1$ and $Z_2$ be real algebraic
  surfaces in $X_{\Delta_k}$ with the Newton polytope $\Delta_k$.  Assume
  that $Z_1$ and $Z_2$ define a complete intersection $Y_k$ in
  $X_{\Delta_k}$.  Then $Y_k$ is not maximal.
\end{proposition}

The proof of Proposition~\ref{ciconterex} relies on the
estimation of the Betti numbers of the
complex and real parts of
the complete intersection $Y_k$ of two
surfaces  whose Newton  polytopes
coincide with $\Delta_k$.

\begin{lemma}\label{lemma-Betti}
Let $Y_k$ be the complete intersection of two
surfaces  in   $X_{\Delta_k}$
whose Newton  polytopes
coincide with $\Delta_k$.
Then
$b_*(Y_k;\CC) = 2k$.
\end{lemma}

\proof
By Corollary~\ref{CorKho}, we
have
$$b_*(Y_k;  \CC)=  2\Vol(\Delta_k) -
\sum_{\Gamma \in \FF_2(\Delta_k)}
\Vol(\Gamma) + 4.$$
So, we get $b_* (Y_k; \CC) = 2k$.
\CQFD

{\bf
Proof of
Proposition~\ref{ciconterex}. }
According to Lemma~\ref{lemma-Betti},
we have $b_*(Y_k; \CC) = 2k$.
Thus, $b_*(Y_k;\ZZ_2) \ge 2k$.

Let   $f_1$  and   $f_2$   be  the   polynomials   defining  the   two
surfaces. Then,
$$f_l(x,y,z) = a_l x + b_l y +c_l z^k + d_l \; (l = 1,2)$$
for some $(a_l, b_l, c_l, d_l)$ in $ \RR^4$.
The change of
variables
$\Lambda_k  : x \mapsto x$, $\Lambda_k  : y \mapsto y$,
$\Lambda_k  : z  \mapsto z^{\frac{1}{k}}$
is a
diffeomorphism of
the first
octant $(\RR_+^*)^3$, where $\RR_+^* = \{x \in \RR : x > 0\}$.
Let $Q_i$ be another
octant, and
$\phi_i$ be the diffeomorphism from $Q_i$ to $(\RR_+^*)^3$
defined by $\phi_i(x,y,z)=(|x|,|y|,|z|)$.
Then $\psi_i =\phi_i^{-1} \circ \Lambda_k \circ \phi_i$
is a diffeomorphism from $Q_i$ to itself. The
diffeomorphism  $\psi_i$ maps  the  zeros  of $f_l$  to  the zeroes  of
${\psi_i}_*(f_l)$ and ${\psi_i}_* (f_l)(x,y,z)= a_l x + b_l y +c_l z +
d_l$. Thus, in each octant, $\RR Y_k$ is diffeomorphic to the intersection
of  two plans. Hence,
the number of connected components of $Y_k$ is
at most $4$. So, $\RR Y_k$ is not maximal for
$k \geq 5$.
\CQFD

The example above should be compared with the following result
in dimension~$2$ which is probably well known but that I couldn't
find in the literature.

\begin{proposition}
Let  $\Delta$ be  a  two-dimensional
polygon. For any positive integers
$\lambda_1$ and $\lambda_2$ there exist algebraic curves
$C_1$ et $C_2$ in $X_\Delta$ such that
\begin{itemize}
\item
the Newton
polygons of $C_1$ et $C_2$
are
$\lambda_1 \cdot \Delta$ and $\lambda_2 \cdot \Delta$,
respectively,
\item
the curves $C_1$ et $C_2$ define a $0$-dimensional
maximal complete intersection
in $X_\Delta$.
\end{itemize}
\end{proposition}

\proof We use here the Cayley trick.  Take any primitive convex
triangulation $\tau$ of $\Delta$. By homothety, $\tau$ induces a
triangulation $\tau_i$ on $\lambda_i \cdot \Delta$.  Put $\Delta_i =
\lambda_i \cdot \Delta$.  Consider the following subdivision
$\delta_0$ of the Cayley polytope $C(\Delta_1, \Delta_2)$.  In the
faces $\hat{\Delta}_1$ and $\hat{\Delta}_2$ of $C(\Delta_1, \Delta_2)$
corresponding to $\Delta_1$ and $\Delta_2$ take the triangulations
$\tau_1$ and $\tau_2$, respectively.  Each $3$-dimensional polytope of
the subdivision $\delta_0$ is the convex hull of a triangle of
$\tau_1$ and a triangle of $\tau_2$ which are the multiples of the
same triangle of $\tau$.  Since $\tau$ is convex, $\delta_0$ is also
convex. Let $\nu_0$ be a convex function certifying the convexity of
$\delta_0$, and let $\nu_1$ be the convex function defined by
$\nu_1(0,1,x,y)= C_1 y + C_2 x$ with $C_1 > C_2 > 0$ and
$\nu_1(1,0,x,y)=0$.  Put $\nu_3 = \nu_1 +\nu_2$.  If $C_1$ is
sufficiently small, the function $\nu_3$ induces the following
refinement $\delta_1$ of $\delta_0$.  Each $3$-dimensional polytope of
$\delta_0$ is subdivided into two cones whose bases are triangles in
$\hat\Delta_1$ and $\hat\Delta_2$, respectively, and a join $J$ of two
edges: one in $\hat\Delta_1$ and the other one in $\hat\Delta_2$.
Take any convex primitive triangulations $\tau_1^\prime$ and
$\tau_2^\prime$ refining $\tau_1$ and $\tau_2$, respectively.  They
define a convex primitive refinement $\delta_2$ of $\delta_1$. Choose
a sign distribution at the vertices of $\delta_2$ and apply the
procedure of the combinatorial patchworking.  Let $J$ be a join of the
decomposition $\delta_1$ described above.  It is triangulated into
primitive tetrahedra $t_i$ and has lattice volume
$\lambda_1\lambda_2$. Each $t_i$ has a symmetric copy containing a
point of the $T$-complete intersection constructed.  Thus, the number
of intersection points obtained is $\lambda_1 \lambda_2 \Vol(\Delta)$
and the complete intersection constructed is maximal.\CQFD

\subsection{Proof of Proposition~\ref{cigencontrex}}

Consider
the
simplex $\sigma_k$ in $\RR^d$
with the vertices
$$\displaylines{(0, 0, \ldots , 0, 0),\> (1, 0, \ldots , 0, 0),\>
(0, 1, \ldots , 0, 0),\>
\ldots ,
(0, 0, \ldots , 1, 0),\cr
\text{\rm and } (1, 1, \ldots , 1, k).}$$
Let $Y_k$ be
a complete intersection  of hypersurfaces  in
$X_{\sigma_k}$ such that all these hypersurfaces
have
 the Newton
polytope $\sigma_k$.
Proposition~\ref{DKcompass}
implies
that
$b_*(Y_k;\ZZ_2)$ tends to infinity when $k$
tends to infinity.

Let $f_1, \ldots , f_n$ be the polynomials defining the hypersurfaces.
Then,
$$f_l(x,y,z) = a_{l,0} + \sum_{i=1}^{d-1} a_{l,i} x_i + a_{l,d}
{x_d}^k \; (l = 1,\cdots , n)$$
for some $(a_{l,0}, \cdots, a_{l,d})$
in $\RR^{d+1}$.  The change of variables $\Lambda_k : x_i \mapsto x_i$
for $i \neq d$, $\Lambda_k : x_d \mapsto {x_d}^{\frac{1}{k}}$ is a
diffeomorphism of the first orthant $(\RR_+^*)^d$.  Let $Q_j$ be another orthant, and $\phi_j$ be
the diffeomorphism from $Q_j$ to $(\RR_+^*)^d$ defined by $\phi_j(x_1,
\cdots , x_d) = (|x_1|, \cdots, |x_d|)$.  Then $\psi_j =\phi_j^{-1}
\circ \Lambda_k \circ \phi_j$ is a diffeomorphism from $Q_j$ to
itself. The diffeomorphism $\psi_j$ maps the zeros of $f_l$ to the
zeroes of ${\psi_j}_*(f_l)$ and ${\psi_j}_* (f_l)(x_1, \ldots , x_d) =
a_{l,0} + \sum_{i=1}^{d} a_{l,i} x_i$.  Thus, in each orthant, $Y_k$
is diffeomorphic to the intersection of $n$ hyperplanes.  Hence,
$b_*(\RR Y_i;\ZZ_2)$ is bounded.

So, there exists  a number $k_0$ such that  for any  $k \ge k_0$
and any complete intersection $Y_k$ in $X_{\sigma_k}$
one has $b_*(\RR Y_k;\ZZ_2) <
b_*(Y_k;\ZZ_2)$.\CQFD

\section{Asymptotically      maximal     families      of     complete
intersections}\label{compassmax}

\subsection{Itenberg-Viro asymptotical statement}\label{IVcompl-int}

The proof of Theorem~\ref{compassmaxth} is based on the following result
of Itenberg and Viro.

\begin{theorem}[I. Itenberg and O. Viro]\label{IVcomp}
Let $\Delta$ be a primitive
$d$-dimensional simplex.
For any $k$-tuple $\lambda_1, \cdots ,\lambda_k$
of natural numbers, there exist piecewise-linear convex functions
$\mu_1 , \ldots , \mu_k$ on $\lambda_1 \cdot \Delta, \ldots ,
\lambda_k \cdot \Delta$, respectively,
and sign distributions at the  vertices of the
corresponding
triangulations of $\lambda_1 \cdot \Delta, \ldots ,
\lambda_k \cdot \Delta$
such that the
real complete intersection in $X_\Delta = \CC P^d$
obtained
via Sturmfels' Theorem~\ref{St} from these data is maximal.
\end{theorem}

In fact, as in Section~\ref{hypassmax}, we use only the following
asymptotical version of Theorem~\ref{IVcomp}.

\begin{theorem}[I. Itenberg and O. Viro]\label{IVtheorem2}
For any positive integers $k$, $m_1, \ldots , m_k$ and $d$ such that
$k \leq d$ and $m_j \geq d + 1$ ($j = 1, \ldots , k$),
there exists a complete intersection~$X$
of multi-degree $(m_1, \ldots , m_k)$ in $\RR P^d$ such that

$$b_*(\R X;\ZZ_2) \geq \sum_{i_1 + \ldots + i_k = d}
\left(\prod_{j=1}^k(m_j - 2)(m_j - 3)
\ldots (m_j - i_j - 1)\right)$$
(the summation is over all possible decompositions
$i_1 + \ldots + i_k = d$ of~$d$ in a sum of $k$ positive integer
numbers).
\end{theorem}

\vskip10pt

The
proof of this asymptotical version
is much simpler than the proof of Theorem~\ref{IVcomp}.
It can be extracted from
\cite{IteVir} and was communicated to us by the authors
of~\cite{IteVir}.
We reproduce their proof here for the completeness.

{\bf Proof of Theorem~\ref{IVtheorem2}.}

The notations used here are those of  Subsection~\ref{IVsection1}.
Take the standard simplices $T^d_{m_1}, \ldots , T^d_{m_k}$
and triangulate the Cayley polytope $C(T^d_{m_1}, \ldots , T^d_{m_k})$
(see Subsection~\ref{Cayley}) in the following way. Let $i_1$, $\ldots$, $i_k$ be
nonnegative integers such that $i_1 + \ldots + i_k = d$,
and put $i_0 = 0$.
For any $j = 1, \ldots , k$ consider the face of $T^d_{m_j}$
with the vertices having the numbers
$$i_1 + \ldots + i_{j-1} + 1, \ldots , i_1 + \ldots + i_j + 1.$$
Denote by $J_{i_1,\ldots, i_k}$ the join of the corresponding faces
of $C(T^d_{m_1}, \ldots , T^d_{m_k})$. The simplices
$J_{i_1,\ldots, i_k}$ (for all the possible choices of nonnegative
integers such that $i_1 + \ldots + i_k = d$) form
a triangulation~$\tau'$ of $C(T^d_{m_1}, \ldots , T^d_{m_k})$.

Take for each simplex $T^d_{m_j}$
the triangulation and the distribution of signs described
in Subsection~\ref{IVsection1}. For the simplices
$\hat T^d_{m_1}, \ldots , \hat T^d_{m_k}$ take
the corresponding triangulations and distributions of signs.
The triangulations
of $\hat T^d_{m_1}, \ldots , \hat T^d_{m_k}$ induce a refinement~$\tau$
of~$\tau'$. Notice that $\tau$ is a primitive triangulation
of $C(T^d_{m_1}, \ldots , T^d_{m_k})$.

\begin{lemma}\label{lemma2} For the complete intersection~$X$
of multi-degree $m_1, \ldots , m_k$
in $\R P^d$ provided according to Proposition~\ref{proposition}
by the
triangulation $\tau$ and the distribution of signs defined above,
one has
$$b_*(\R X;\ZZ_2) \geq \sum_{i_1 + \ldots + i_k = d}
\left(\prod_{j=1}^k(m_j - 2)(m_j - 3)
\ldots (m_j - i_j - 1)\right)$$
(the summation is over all the possible decompositions
$i_1 + \ldots + i_k = d$ of~$d$ in a sum of $k$ positive integer
numbers).
\end{lemma}

{\bf Proof. }
We define a collection of narrow cycles $c_i$, $i \in I$
of $\widetilde H$.
The families of narrow cycles of $\widetilde H$ are indexed by
the decompositions
$i_1 + \ldots + i_k = d$ of~$d$ in a sum of $k$ positive integer
numbers.

Fix a decomposition ${\cal I}: i_1 + \ldots + i_k = d$ of~$d$,
where $i_1$, $\ldots$, $i_k$ are positive integers.
The initial data for constructing a narrow cycle
of the corresponding family consist
of narrow cycles
$c_{(j)} \subset \widetilde H^{\cal I}_j$, $j = 1, \ldots , k$,
constructed in Subsection~\ref{IVsection1}
for the hypersurface $\widetilde H^{\cal I}_j$
in $\widetilde T^{i_j}_{m_j}$
produced via the combinatorial patchworking
by the triangulation and distribution of signs described
in Subsection~\ref{IVsection1}.

The $i_j$-dimensional face $\Delta^{i_j}$ of $T^d_{m_j}$
with the vertices having the numbers
$$i_1 + \ldots + i_{j-1} + 1, \ldots , i_1 + \ldots + i_j + 1$$
are naturally
identified with $T^{i_j}_{m_j}$ via the linear
map ${\cal L}^{i_j}: \Delta^{i_j} \to T^{i_j}_{m_j}$
sending
the vertex with number $i_1 + \ldots + i_{j-1} + r$ of $\Delta^{i_j}$
to the vertex with number~$r$ of $T^{i_j}_{m_j}$.
The map ${\cal L}^{i_j}$ is simplicial
with respect to the chosen triangulations of $\Delta^{i_j}$
and $T^{i_j}_{m_j}$.
Denote by $\Delta^{i_j}_*$
the union of the symmetric copies of $\Delta^{i_j}$
under the reflections
with respect to coordinate hyperplanes
$\{x_i = 0\}$ in $\R^d$, where
$i = i_1 + \ldots + i_{j-1} + 2, \ldots, i_1 + \ldots + i_j + 1$,
and compositions
of these reflections.
The natural extension of ${\cal L}^{i_j}$ to $\Delta^{i_j}_*$
identifies $\Delta^{i_j}_*$ with $(T^{i_j}_{m_j})_*$
and respects the chosen triangulations.
We also denote this extension by ${\cal L}^{i_j}$.
Denote by $\hat \Delta^{i_j}_*$ the union of faces
of $\hat T^d_{m_j}$ corresponding to $\Delta^{i_j}_*$,
and by $\hat{\cal L}^{i_j}$ the corresponding
map from $\hat\Delta^{i_j}_*$ to $(T^{i_j}_{m_j})_*$.
Put $\hat c_{(j)} = (\hat{\cal L}^{i_j})^{-1}(c_{(j)})$.

Let $b_{(j)} \subset \widetilde T^{i_j}_{m_j} \setminus
\widetilde H^{\cal I}_j$
be the dual cycle of~$c_{(j)}$.
Put $\hat b_{(j)} = (\hat{\cal L}^{i_j})^{-1}(b_{(j)})$.
Consider the symmetric copies of $\hat b_{(1)}, \ldots , \hat b_{(k)}$
under the reflections
with respect to coordinate hyperplanes
$\{x_i = 0\}$ in $\R^{k+d}$ where
$i = k + 1, \ldots, k + d$,
and compositions
of these reflections.
Among these symmetric copies
there exist copies $\hat b'_{(1)}, \ldots , \hat b'_{(k)}$
of $\hat b_{(1)}, \ldots , \hat b_{(k)}$, respectively, such that
\begin{itemize}
\item the join $\hat b'_{(1)}* \ldots * \hat b'_{(k)}$
is the union of simplices
of $\tau_*$,
\item all the vertices of $\hat b'_{(1)}* \ldots * \hat b'_{(k)}$ have
the same sign.
\end{itemize}
Let $\hat c'_{(1)}, \ldots , \hat c'_{(k)}$
be the corresponding symmetric
copies of $\hat c_{(1)}, \ldots , \hat c_{(k)}$,
respectively.
Then, take the intersection
$B \cap (\hat c'_{(1)}* \ldots * \hat c'_{(k)})$ as a narrow cycle
of $\widetilde H$.

The number of narrow cycles in the family indexed by~${\cal I}$
is at least
$$\prod_{j=1}^k(m_j - 2)(m_j - 3)
\ldots (m_j - i_j - 1).$$
Thus, the total number of constructed narrow cycles
in $\widetilde H$
is at least
$$\sum_{i_1 + \ldots + i_k = n}
\left(\prod_{j=1}^k(m_j - 2)(m_j - 3)
\ldots (m_j - i_j - 1)\right)$$
(the summation is over all the possible decompositions
$i_1 + \ldots + i_k = d$ of~$d$ in a sum of $k$ positive integer
numbers).
The linear independence of the narrow cycles of
a hypersurface $H^l_{m_j} \subset T^l_{m_j}$ for any
$1 \leq l \leq d$ and any $1 \leq j \leq k$ implies
the linear independence of the narrow cycles constructed
in $\widetilde H$.
\CQFD

\begin{remark}\label{rmqsimplex}
Denote by
$Y^\sigma_{m_1, \ldots , m_k}$
the complete intersection
constructed in Lem\-ma~\ref{lemma2}.
Then, the family $\{Y^\sigma_{m_1, \ldots , m_k}\}_{m_1, \ldots , m_k}$
is asymptotically maximal.
\end{remark}

\proof
Note that
$$\sum_{i_1 + \ldots + i_k = d}\left(\prod_{j=1}^k(m_j - 2)
(m_j - 3) \ldots (m_j - i_j - 1)\right)$$
is equivalent to the mixed volume of
$T^d_{m_1}, \ldots ,  T^d_{m_k}$. Thus,
by Proposition~\ref{DKcompass},
$b_*(\RR Y_{m_1, \ldots , m_k};\ZZ_2)$
is equivalent to $b_*(Y_{m_1, \ldots ,  m_k};\ZZ_2)$,
when all $m_i$'s tend to infinity.

\subsection{Proof of Theorem~\ref{compassmaxth}}

Let $\tau$ be a primitive convex triangulation of $l \cdot \Delta$,
and $(\lambda_1, \cdots , \lambda_k)$ be a $k$-tuple of positive
integers.  Denote by $\Delta_{\lambda_i}$ the polytopes $\lambda_i l
\cdot \Delta$.  We can assume that $\lambda_i$ is greater than $d+1$
for any $i$.

Let $\delta$ be a $d$-dimensional simplex of the triangulation $\tau$.
Denote by $\hat\delta_1$, $\ldots$, $\hat\delta_k$ the corresponding
simplices in $\hat\Delta_{\lambda_1}$, $\ldots$,
$\hat\Delta_{\lambda_k}$, respectively.  Subdivide the Cayley polytope
$C(\Delta_{\lambda_1}, \ldots , \Delta_{\lambda_k})$ into convex hulls
of $\hat\delta_1$, $\ldots$, $\hat\delta_k$, where $\delta$ runs over
all $d$-dimensional simplices of $\tau$.  For a $d$-dimensional
simplex $\delta$ of $\tau$, put $\delta_i=\lambda_i \cdot \delta$ and
$\delta_i^\prime=(\lambda_i-(d+1)) \cdot \delta$, where $i = 1, \ldots
, k$.

For any
$d$-dimensional simplex $\delta$ of
$\tau$, take the triangulation of
$C(\delta_1^\prime, \ldots , \delta_k^\prime)$
and the distribution of signs at the vertices of this triangulation
described in the proof of Theorem~\ref{IVcompl-int}.
Extend the triangulations of the Cayley polytopes
$C(\delta_1^\prime, \ldots , \delta_k^\prime)$ to
a primitive convex triangulation $\hat\tau$ of
$C(\Delta_{\lambda_1}, \ldots , \Delta_{\lambda_k})$ in
the same way as it was done in Subsection~\ref{assmax-section}.
Extend also the distributions of signs at the integer
points
of polytopes $C(\delta_1^\prime, \ldots , \delta_k^\prime)$
to some distribution
of signs $\hat D$ at the vertices of $\hat\tau$.

Let $Y_{\lambda_1, \ldots , \lambda_k}$ be
the complete intersection in $X_\Delta$ obtained
via Theorem~\ref{proposition}
from $\hat\tau$ and $\hat D$.

\begin{proposition}
The  family of complete intersections
$Y_{\lambda_1, \ldots , \lambda_k}$
constructed above is asymptotically maximal.
\end{proposition}

\proof By the construction, we have $b_*(\RR Y_{\lambda_1, \ldots ,
  \lambda_k};\ZZ_2) \ge \Vol(l \cdot\Delta) \cdot n_{\lambda_1, \ldots ,
  \lambda_k}$, where $n_{\lambda_1, \ldots , \lambda_k}$ is the number
of narrow cycles in each $C(\delta_1^\prime, \ldots ,
\delta_k^\prime)$.  Note that $n_{\lambda_1, \ldots , \lambda_k}$ is
equivalent to $b_*(\RR Y^\sigma_{\lambda_1, \ldots ,
  \lambda_k};\ZZ_2)$ when all numbers $\lambda_1, \ldots , \lambda_k$
tend to infinity. So, by Proposition~\ref{DKcompass} and
Remark~\ref{rmqsimplex}, we obtain that $b_*(\RR Y_{\lambda_1, \dots ,
  \lambda_k};\ZZ_2)$ is equivalent to $b_*( Y_{\lambda_1, \dots ,
  \lambda_k};\ZZ_2)$ when the numbers $\lambda_1, \ldots , \lambda_k$
tend to infinity.  \CQFD

\providecommand{\bysame}{\leavevmode\hbox to3em{\hrulefill}\thinspace}
\providecommand{\MR}{\relax\ifhmode\unskip\space\fi MR }
\providecommand{\MRhref}[2]{%
  \href{http://www.ams.org/mathscinet-getitem?mr=#1}{#2}
}
\providecommand{\href}[2]{#2}

\end{document}